%title: critical conditions for thermal explosion
%author: barenblatt, chorin, kast
%date: 9/10/97

\input amstex
\input amsppt.sty
\magnification=\magstep1
\hsize = 6.5 truein
\vsize = 9 truein
\NoBlackBoxes
\TagsAsMath
\NoRunningHeads

%The next three lines are so the logo doesn't print.
\catcode`\@=11
\redefine\logo@{}
\catcode`\@=13

%Section macros
\newskip\sectionskipamount
\sectionskipamount = 24pt plus 8pt minus 8pt
\def\sectionskip{\vskip\sectionskipamount}
\define\sectionbreak{%
    \par  \ifdim\lastskip<\sectionskipamount
    \removelastskip  \penalty-2000  \sectionskip  \fi}
\define\section#1{%
    \sectionbreak   %Encourages a page break, else inserts 24pt glue
    \subheading{#1}%
    \bigskip
    }

\font\bigbf=cmbx10 scaled\magstep 1

\addto\tenpoint{\normalbaselineskip 18truept \normalbaselines}
\def\label#1{\par%
    \hangafter 1%
    \hangindent 1.25 in%
    \noindent%
    \hbox to 1.25 in{#1\hfill}%
    \ignorespaces%
    }

%Operator name macros
\define\op#1{\operatorname{\fam=0\tenrm{#1}}} %for text in math mode use
                                              %\op{...}

%QED box, from the TeXbook, p. 106.

    \let    \< = \langle
    \let    \> = \rangle
    \define      \a     {\alpha}
    \redefine   \b      {\beta}
 \redefine \dt {\cdot}
    \redefine   \d      {\delta}
    \redefine   \D      {\Delta}
    \define     \e      {\varepsilon}
    \define      \g     {\gamma}
    \define     \G      {\Gamma}
 \redefine \k  {\kappa}
    \redefine   \l      {\lambda}
    \redefine   \L      {\Lambda}
    \redefine   \ni     {\noindent}
    \redefine   \var    {\varphi}
    \define     \s       {\sigma}
    \define     \th     {\theta}
    \redefine   \O      {\Omega}
    \redefine   \o      {\omega}
    \define      \z     {\zeta}
    \redefine   \i      {\infty}
    \define      \p     {\partial}

\document
\baselineskip=16 pt

\centerline{\bigbf The Influence of the Flow of the Reacting Gas
}
\centerline{{\bigbf on the Conditions for a Thermal Explosion
}\footnote{Supported in part by the Applied Mathematical
Sciences subprogram of the Office of Energy Research, U.S.
Department of Energy, under contract DE--AC03--76--SF00098, and in
part by the National Science Foundation under grant DMS94--14631
.}}
\bigskip\bigskip
\centerline{\smc G.~I.~Barenblatt, A. J.~Chorin, and A.~Kast}

\medskip
\centerline{Department of Mathematics}
\centerline{Lawrence Berkeley National Laboratory}
\centerline{Berkeley, California 94720--3840} \medskip

\bigskip\bigskip
{\bf Abstract.} The classical problem of thermal
explosion is modified so that the chemically active
gas is not at rest but is flowing in a long cylindrical pipe. Up to a
certain section the
heat-conducting walls of the pipe are held at low temperature
so that the reaction rate is small and there is no heat release; at that
section the ambient temperature is increased and
an exothermic reaction begins. The question is whether a slow
reaction regime will be established or a thermal explosion will
occur.

The mathematical formulation of the problem is presented.
It is shown that when the pipe radius is larger than a
critical value,
the
solution of the  new problem exists only up to a certain
distance along the axis.
The critical radius is determined by conditions in a problem with
a uniform axial temperature.
The loss of existence is
interpreted as a thermal explosion; the critical
distance is the safe reactor's length. Both
laminar and developed turbulent flow regimes are considered. In a
computational experiment
the loss of the existence appears as a
divergence of a numerical procedure;
numerical calculations
reveal asymptotic scaling laws with simple powers for the critical
distance.

\newpage
\ni{\bf 1. \ Introduction}

The classical problem of a thermal explosion in a long cylindrical
vessel containing a reacting gas at rest was first formulated and solved
by D.~A.~Frank-Kamenetsky in 1939 (see (${}^1$)). The basic result is
that there is no spontaneous thermal explosion if the radius of the vessel
$r_0$ is less than a critical value
$r_*$ determined by the gas properties and the properties of the
exothermal chemical reaction.
If $r_0 > r_*$, a quiet evolution  of the chemical
reaction is impossible and a thermal explosion occurs.

In the present paper the following question is asked: What will
happen if the long cylindrical vessel is in fact a pipe where the
chemically active gas is flowing? At a certain section of the pipe
(see Figure 1) we suddenly increase the temperature of the walls
and the rate of the exothermic reaction increases. How will the
critical conditions for thermal explosion  be affected by the gas flow?
In the present paper we provide a mathematical formulation of the problem.
Our first result is rather obvious: if the radius of the pipe is
less than critical at the increased wall temperature, nothing will
happen; the reactor will be safe for arbitrary length. The second
result seems to be deeper: If the radius of the pipe is larger
than critical, there exists a distance along the
pipe up to which the solution of the mathematical
problem exists, and after which it ceases to exist. This critical length
depends on the flow conditions and
determines the ``safe" length of the reactor, up to which the
thermal explosion  does not happen.

We want to mention here the remarkable paper (${}^2$) where the
influence of a flow on the critical conditions for thermal explosion  was
first
investigated. The flow under consideration in (${}^2$) was
different--free convection-- and therefore the mathematical problem
was also different from ours.

\bigskip\bigskip\ni{\bf 2. Problem formulation: Laminar flow}

We make the same basic assumptions concerning the exothermic
chemical reaction as in the classical theory of thermal explosion  (see
(${}^1$)): Arrhenius dependence of the reaction rate on the temperature, large
activation energy, and neglect of the longitudinal heat transfer in
comparison with the contribution of the advection. Then the
equation of the energy balance takes the form:
$$
cv(r)\p_z T=\frac 1r \p_r (rk\p_r T)+Q\s(T)e^{-E/RT} \ . \tag 1
$$
Here $T$ is the absolute temperature of the gas, $c$ is the heat capacity,
$z$ is the longitudinal coordinate reckoned from the section where the
temperature was elevated, $v(r)$ is the longitudinal fluid velocity,
$r$ is the current radius, $k$ is the heat conductivity corresponding to
the flow regime, $Q$ is the thermal effect of the reaction, $\s(T)$ is
the pre-exponential factor in the Arrhenius reaction rate as a function of the
temperature- a slow function of the temperature, $E$ is the
activation energy, and $R$ is the universal gas constant. Applying the
Frank-Kamenetsky high activation energy approximation (see
(${}^1$)), and going to dimensionless variables, we reduce
equation (1) to the form
$$
V(\rho)\p_\zeta u=\frac{1}{\rho} \p_{\rho}\rho
\left( \frac{\k}{\k_0} \p_{\rho}u\right) +\l^2 e^u \ . \tag 2
$$
Here $\rho=r/r_0$ is the dimensionless radius, $\k$ the temperature
diffusivity appropriate for the flow regime under consideration,
$\k_0$ the molecular temperature diffusivity, $V(\rho)=v(r)/\bar
v$, \ $\bar v=G/\pi r^2_0$, is the average velocity, $G$ the flow
discharge rate, $\zeta=z/(G/\pi\nu_0)Pr$ the dimensionless
longitudinal coordinate, $\nu_0$ the molecular kinematic viscosity,
$Pr=\nu_0/\k_0$ the Prandtl number, a characteristic gas constant.
Furthermore, $\l=r_0/\ell$, where
$\ell=(e^{\phi}\kappa_0T_0c/Q\phi\s(T_0))^{\frac 12}$, \ $\phi=E/RT_0$,\
$u=E(T-T_0)/RT^2_0$, $T_0$ being the wall temperature at $z > 0$.
For a gas temperature uniformly distributed along the $z$-axis
in a laminar flow
$\p_zT=0$, $\k=\k_0$ and we return to the classical problem
formulation (see (${}^1$)):
$$
\frac{1}{\rho} \ \frac{d}{d\rho}\left(\rho \frac{du}{d\rho}\right)
+\l^2e^u=0 \ , \qquad
u(1)=0, \qquad
\lim_{\rho\to 0} \left(\rho \frac{du}{d\rho}\right)=0. \tag 3
$$

In the case of laminar motion it is easy to show that
$V(\rho)=2(1-\rho^2)$, therefore equation (2)
assumes the form
$$
(1-\rho^2)\p_\zeta u=\frac{1}{\rho} \p_\rho \rho\p_\rho u+\l^2e^u
\tag 4
$$
(we redefined $\zeta: \zeta=z/(2G/\pi\nu_0)Pr)$. The boundary
conditions for the solution $u(\rho,\zeta)$ take the form
$$
u(\rho,0)=0, \qquad
u(1,\zeta)=0, \qquad
\lim_{\rho\to 0} (\rho\p_\rho u(\rho,\zeta))=0. \tag 5
$$
The first condition (5) reflects the assumption that the temperature of
the flowing gas in the whole inlet section $z=0$ is equal to the
temperature of the wall for $z>0$. This is not essential and can be
replaced by an
arbitrary inlet condition $u(\rho,0)=u_0(\rho)$. The second and
third conditions are the same as in the classical problem
formulation. It is clear that $\p_{\zeta}u$ is positive,
therefore for $\l < \l_{\text{cr}}=\sqrt{2}$, the solution of the
problem does exist for  $0 \leq\zeta <\i$. For  $\l >
\l_{\text{cr}}$ the solution to the problem (4)--(5) does exist
only at  $0 \leq\zeta\leq\zeta_0$, where $\zeta_0$ is a certain
finite positive constant which depends on $\l$. The critical
dimensional length is
$$
z_0 =r_0 Re \ Pr \ \zeta_0, \qquad
Re  =\frac{\bar vd}{\nu_0}=\frac{2\bar vr_0}{\nu}.
\tag6$$ In
section 4 the numerical solution of the problem (4)--(5) will be presented.

\bigskip\bigskip\ni{\bf 3. Problem formulation:
Developed turbulent flow}

In this case we use the power law proposed earlier
(see (${}^3$)), according to which the longitudinal velocity
distribution takes the form
$$
v(r)=v_* C(\a)\left( \frac{v_*y}{\nu_0}\right)^{\a} \ , \qquad
\a=\frac{3}{2\ln Re} \ , \qquad
C(\a)=\frac{\sqrt{3} +5\a}{2\a} \ .\tag 7
$$
Here $y=r_0-r$ is the distance from the wall, $Re=2\bar v
r_0/\nu_0$ is the Reynolds number, and $v_*$ is the ``friction" velocity.

Equation (2) in this case is reduced to the form
$$
(1-\rho)^\a \p_\zeta u=\frac{1}{\rho} \p_\rho \rho^2
(1-\rho)^{1-\a} \p_\rho u+\l^2e^u \ . \tag 8
$$
The first two boundary conditions (5) preserve their form in
the turbulent case, however the third condition --- the condition
of zero influx at the axis --- takes the form
$$
\lim_{\rho\to 0} [\rho^2(1-\rho)^{1-\a}\p_\rho u]=0. \tag 9
$$
To obtain (8),(9) an expression for the turbulent temperature diffusivity
$\k$ is needed. This quantity is determined by a relation
$\k=A(Re,Pr)\nu_T$, where  $\nu_T$ is the turbulent kinematic
viscosity and $A(Re,Pr)$ is a dimensionless factor which we do not
necessarily assume
to be equal to unity. By definition
$\nu_T=v^2_*\rho/\p_y v$, so that
$$
\k=\k_0\frac{A(Re,Pr)Pr}{\alpha
C(\alpha)2^{1-\alpha}}\left(\frac{v_*d}{\nu_0}\right)^{1-\a}
\rho(1-\rho)^{1-\alpha}.
\tag 10
$$
The quantity $(v_*d/\nu)^{1-\alpha}$ is a function of $\alpha$, or, which
is the same, of
the Reynolds number $Re$ (see (${}^3$)). Substituting this expression into
(2), and redefining $\zeta$ and $\ell$ by
$$
\aligned
\zeta& =\frac{z\cdot 2^{2\alpha}\cdot
A(Re,Pr)}{r_0\alpha C^2(\alpha)}
\left(\frac{v_*d}{\nu_0}\right)^{-2\alpha}, \\
\ell^2 & =\left(\frac{e^{\phi}\kappa_0T_0c}{Q\phi\sigma(T_0)}
\right)\frac{PrA(Re,Pr)}{2^{1-\alpha}\alpha
C(\alpha)} \left(\frac{v_*d}{\nu_0}\right)^{1-\alpha},
\endaligned   \tag 11$$
we obtain the basic equation (8). Condition (9) is obtained from (10) and
the from the
condition of zero flux at the axis $\rho=0.$

The critical dimensional length $z_0$ is therefore determined by the relation
$$z_0=r_0\alpha
C^2(\alpha)\left(\frac{v_*d}{\nu_0}\right)^{2\alpha}
\frac{\zeta_0}{2^\alpha A(Re,Pr)}.\tag12$$

To demonstrate the difference between condition (9) and the
corresponding condition for the laminar case, consider the equation
$$
\frac{1}{\rho} \ \frac{d}{d\rho} \rho^2(1-\rho)^{1-\a}
\frac{du}{d\rho}+\l^2e^u=0, \tag 13
$$
which, in the subcritical case, describes the temperature distribution at
large distances from the inlet section, where the longitudinal temperature
variation is negligible. If condition (9) holds, $u(0)$ is
finite, so that integrating (13) we obtain
$$
\rho^2(1-\rho)^{1-\a} \frac{du}{d\rho} =-\l^2\int^{\rho}_0
\rho e^{u(\rho)}d\rho; \qquad
\left(\frac{du}{d\rho}\right)_{\rho=0}=-\frac{\l^2}{2} \ e^{u(0)} \
. \tag 14
$$
Thus the integral curves of (13) that satisfy condition (9) form a
one-parameter family having $u(0)$ as parameter. This family
possesses an envelope, and the value of the parameter $\l$ which
corresponds to the intersection of the envelope with the
$\rho$-axis at $\rho=1$ is the critical value.

In the following section the value of $\l_{\text{cr}}$ will be
determined numerically, in the case of turbulent flow, as a function of the
Reynolds number, and for values $\l > \l_{\text{cr}}$ the
critical dimensionless length $\zeta_0$ will also be determined.

\bigskip\bigskip\ni{\bf 4. Numerical solutions}

It is convenient to use the coordinates $\rho, \xi=\zeta/\zeta_0$ in the
numerical solution of the problems
(4)--(5) and (8),(5),(9)
so that the domain of
integration is fixed: $0\leq\rho\leq 1$, $0\leq\xi\leq 1$;
in these coordinates
$\zeta_0$ appears as a parameter in the differential equations.
The resulting equations are discretized and yield a system of
algebraic equations for $u$ at grid points. This system is
solved by Newton's method. If Newton's method converges for a
particular pair $\l,\zeta_0$ we decide that the pair is inside the
existence region in the $\l,\zeta_0$-plane. If the method does not
converge, the pair is outside.

Newton's method is iterative and an initial condition for iteration
must be chosen. If this condition is close to an actual solution
then the method is rapidly convergent. If the initial condition is
not close to the actual solution then the method can fail to
converge even when the solution does exist. To pick the initial
conditions, we first assume that for any $\l$, a solution can be
found if the value of parameter $\z_0$ is sufficiently small. We
pick $\l$, compute a solution for some $\z_0$ sufficiently small
so that the
solver converges, and solve again starting from the initial
conditions given by the preceding solution. When $\z_0$ becomes so
large that the solver fails to converge, we say that we have
crossed the existence boundary. In the laminar case, whenever
$0 < \l <\sqrt{2}$, our method identifies convergent solutions for
all values of $\z_0$ tested. When $\l > \sqrt{2}$, our method
identifies convergent solutions for all values of $\z_0$ less
than a critical value $\z_0(\l)$.

Figure 2 shows the graph of $\z_0(\l)$ for the laminar case.
The salient feature is the divergence of $\z_0(\l)$ at the point
$\l=\sqrt{2}$: \ $\z_0(\sqrt{2})=\i$. Similar calculations were
also performed for the turbulent case. Figure 3 shows the graph
of $\l^2_{\text{cr}}$ as a function of the Reynolds number.
Figure 4 shows the graph of the critical values $\zeta_0(\l)$.
A detailed description of the numerical work will be provided elsewhere.

\bigskip\bigskip\ni{\bf 5. Conclusion}

The formulation of a new mathematical problem for chemically
active gas flow in a pipe is presented. It is shown that for a
supercritical regime, where the radius of the pipe is larger than a
critical value, the solution of this problem exists only for an
interval of bounded
length  along the pipe axis. The loss of the existence
of the solution is naturally interpreted as a thermal explosion . The
dimensional
critical length $\z_0$ is determined both for laminar and for developed
turbulent flows.

In the case of laminar flow the results of the numerical
computation of $\z_0(\l)$ ($\l$ is the parameter proportional to the
radius of the pipe) are presented in Figure 2. The asymptotic
branch of the curve  $\z_0(\l)$ at large $\l$ corresponds to the
scaling law  $\z_0(\l)=1.77\l^{-11/4}$. The results of the numerical
computation of $\l_{\text{cr}}$ as a function of the Reynolds
number for developed turbulent flow are presented in Figure 3,
and the function  $\z_0(\l)$ for developed turbulent flow at
various values of the Reynolds number  are presented in Figure 4.
An interesting property of this graph is the collapse of the curves
corresponding to different values of the Reynolds number in a very
large (3 orders of magnitude) range of $Re$ starting from rather
moderate values of $\l$. The asymptotic branch of the curve
 $\z_0(\l)$ at large $\l$ corresponds to the scaling law
 $\z_0(\l)=0.78\l^{-9/4}$. The origin of the observed scaling laws
remains unknown to us.

\bigskip\bigskip
\centerline{\bf References}

(${}^1$) Zeldovich, Ya.B., Barenblatt, G.I., Librovich, V.B., and
Makhviladze, G.M. (1985). The Mathematical Theory of Combustion and
Explosions (Consultants Bureau, New York, London).

(${}^2$) Kagan, L., Berestycki, H., Joulin, G., and Sivashinsky, G.
(1997),
{\it Combust. Theory Modelling} {\bf 1}, 97.

(${}^3$) Barenblatt, G.I., Chorin, A.J., and Prostokishin, V.M.
(1997),
{\it Appl. Mech. Rev.} {\bf 50}, 413.

\newpage

\centerline{\bf Figure Captions}

\ni Figure 1. Chemically active gas flow in a pipe (schematic).

\bigskip\bigskip\ni Figure 2. The dimensionless critical length $\z_0$ as a
function of the basic parameter $\l$ proportional to the pipe
radius, in the case of laminar flow. The critical length tends to
infinity as $\l$ approaches the critical value
$\l_{\text{cr}}=\sqrt{2}$ from above, and exhibits the scaling law
$\z_0\sim\l^{-11/4}$ at large supercritical values of $\l$.

\bigskip\bigskip\ni Figure 3. The critical value of the basic parameter
$\l_{\text{cr}}$ for turbulent flow as a function of the
Reynolds number $Re$.

\bigskip\bigskip\ni Figure 4. The dimensionless critical length $\z_0$ as a
function of the basic parameter $\l$ in the case of
developed turbulent flow. Note the collapse of the curves starting
from rather moderate values of $\l$ in a large range of Reynolds
number, and the scaling law $\z_0\sim\l^{-9/4}$ at large
supercritical values of $\l$.

\enddocument